# An Engineering and Statistical Look at the Collatz (3n + 1) Conjecture

**Brian Mohan Gurbaxani**
**14 Mar 2021**
**(GurbaxaniMath@gmail.com)**

**Abstract**

The famous (3n + 1) or Collatz conjecture has admitted some progress over the last several decades towards the conclusion that the conjecture is true (i.e. that all Collatz sequences will eventually reach a value of one), but has stubbornly resisted a proof. Many professional mathematicians have applied an impressive array of machinery to its analysis, and discovered some limits or boundaries to where violations of the Collatz conjecture must be, if they exist. Here we attempt to look at the conjecture differently than most mathematicians or number theorists, i.e. we will study the conjecture from the points of view of a statistician/data scientist and of an engineer. As a statistician or data scientist, we look at Collatz sequences as if they were sequences found in nature, like a collection of time series produced by some natural process, ignoring for the time being their completely deterministic origins. As an engineer, we try to tinker with Collatz sequences to see what makes them tick, and engineer -- tweak and build – changes to the sequences that likewise do interesting things, after first specifying exactly what we mean by "interesting". Although these analyses do not yield a proof of the Collatz conjecture, like other efforts they present reasons to think that the conjecture is probably true, and the hope is that the analyses of sequences similar to Collatz sequences will tell us something about the nature of Collatz conjecture itself, and/or lead us into interesting and fruitful new directions. Indeed, many of the phenomenon revealed in the sequences of these Collatz-like programs are arguably more interesting than, and suggest conjectures that are as interesting as, the Collatz conjecture itself.

**Disclaimer up front**: Because I am not a professional mathematician, some of the results derived here may be contained elsewhere in the literature without my knowledge. Please forgive me in advance for any failure to cite those prior results – the oversight was not intentional. I would be happy to reference those results in an updated manuscript.

**Introduction**

The Collatz conjecture states simply that if you start with any positive integer n0, then divide it by 2 if n0 is even or multiply by 3 then add 1 if n0 is odd, and keep repeating the same process with each new result, eventually you will always arrive at the number 1. As a convention, we will assume (3n + 1)/2 as a single step for odd numbers in the sequence (since multiplying any odd number by 3 and adding 1 inevitably produces an even number), and refer to them loosely as a "3/2 step increase". Steps applied to the even numbers will be referred to as a "1/2 step". Collatz sequences are known for their highly variable lengths even given very similar n0 (for example n0 = 27 has 71 iterations to get to 1 as we define iterations here, but n0 = 29 gets to 1 in only 13 iterations), and seemingly random strings of 3/2 steps and 1/2 steps. The sequence for n0 = 29 is shown below as an example:

29 -> 44 -> 22 -> 11 -> 17 -> 26 -> 13 -> 20 -> 10 -> 5 -> 8 -> 4 -> 2 -> 1

Some of what has been discovered before about Collatz sequences (Lagarias 1985, Lagarias 1990, Shaw 2006, Tao 2019, OEIS A006577) is summarized below:

1) The truth of the conjecture has been verified for at least all positive n0 < $10^{20}$
2) If a repeating cycle or loop exists other than the one at the bottom of the number line (1 -> 4 -> 2 -> 1), it must have > 250,000 iterations. A simple modification to the Collatz sequence of 3n -1 (subtract 1 instead of add) yields two more cycles or loops which have minima at 5 and 17 for the modified Collatz program, i.e. three loops for 3n-1, one for 3n+1.
3) No sequences are expected to keep increasing to infinity because of the general expectation that sequences decrease on average by a factor of ¾ every 2 steps (assuming 3/2 steps and 1/2 steps are equally likely to occur).
4) There are no multiples of 3 inside the sequences themselves, i.e. n0 containing a factor of 3 are terminal nodes of Collatz sequences (not including their multiples of powers of 2).
5) Sequence lengths (s(n0) for any given starting integer n0) have a high degree of persistence, as nearly half the time s(n0+1) = s(n0).

This is not an exhaustive list of prior findings, but these are some of the findings that are most often referenced in discussions of the Collatz conjecture. We will look at some other properties of Collatz sequences, examining the sequences as a computer engineer might – in terms of the binary operations that create the sequences and in terms of graph theory – before tinkering with the basic sequence generating program to discover other sequences that behave like the Collatz generator, but are arguably even more interesting.

**Results**

The results will be presented as a series of questions and answers as a better way to summarize the findings and (hopefully) so the reader can easily find the questions they are most interested in. The questions and answers are generally, but not necessarily, presented in the chronological order in which they were posed and/or discovered (or re-discovered if that is the case, by the author). Following the questions/answers some follow on queries and conjectures will be addressed/restated/summarized. Note: in what follows, as a convention taken from early programming languages (like FORTRAN), I will to use the letters {i – n} inclusive for integer variables.

**Q1**: A look at the sequence of iterations for any given n0 shows an almost random sequence of 3/2 seps and 1/2 steps, with some sequences containing long sub-strings of consecutive 3/2 steps that create very large n before the sequence eventually comes down to 1. Is there a limit to how many consecutive 3/2 increase steps in a row a Collatz sequence can contain or can they increase without limit?

**A1**: Looking at the sequences in binary, and it is easy to see that sequences starting n0 = $2^k – 1$ are front loaded with k 3/2 steps in a row, thus the initial run-up in Collatz sequences can increase indefinitely as n0 increases. This also makes it easy to approximate the highest n in a given Collatz sequence starting with an n0 of this form.

In Figure 1 below, we see how the Collatz procedure, visualized on binary numbers, inevitably produces a run up in value due to several least significant 1 bits in a row. When the least significant bit is a 1, that indicates an odd integer, and the next step in the sequence involves a 3n + 1 followed by a divide by 2.

A multiply by 3 is equivalent to a multiply by two then adding the result to itself, i.e. 3n = 2n + n. The multiply by 2 is simply a register shift left in binary, and the +1 simply puts a 1 bit in the least significant bits after the shift left. A series of consecutive 1 bits followed by at least one 0 bit will continue to be shortened by a 1 bit with each successive (3n+1)/2 step until the runup is complete. The intervening 0 bits may grow or not but there will always be at least one of them preserved in the process, acting as a type of buffer ensuring that the runup caused by the consecutive 1 bits lasts only for the number of those consecutive bits. For example, Figure 1 shows the Collatz procedure in binary for two consecutive 3/2 steps. Because the starting number has ten consecutive 1 bits at the beginning of its binary sequence, the runup will consist of ten 3/2 steps in a row before an intervening 1/2 step or multiple 1/2 steps. Thus, for n0 = $2^k - 1$, where n0 consists of only 1 bits in binary, the initial runup is k consecutive 3/2 steps. So the ratio of the highest number in a Collatz sequence and its starting point n0 can grow exponentially, depending on how many 1 bits front load the n0 when viewed in binary.

| | | | | | | | | | | | | | | | | | | | | | | | | | | |
|---|---|---|---|---|---|---|---|---|---|---|---|---|---|---|---|---|---|---|---|---|---|---|---|---|---|---|
| | n | : | | | 1 | 0 | 1 | 0 | 1 | 1 | 0 | 0 | 1 | 1 | 1 | 0 | 1 | 1 | 1 | 1 | 1 | 1 | 1 | 1 | 1 | 1 |
| + | 2n+1 | : | | 1 | 0 | 1 | 0 | 1 | 1 | 0 | 0 | 1 | 1 | 1 | 0 | 1 | 1 | 1 | 1 | 1 | 1 | 1 | 1 | 1 | 1 | 1 |
| | 3n+1 | = | 1 | 0 | 0 | 0 | 0 | 0 | 0 | 1 | 1 | 0 | 1 | 1 | 0 | 0 | 1 | 1 | 1 | 1 | 1 | 1 | 1 | 1 | 1 | 0 |
| | /2 | = | | 1 | 0 | 0 | 0 | 0 | 0 | 0 | 1 | 1 | 0 | 1 | 1 | 0 | 0 | 1 | 1 | 1 | 1 | 1 | 1 | 1 | 1 | 1 |
| + | 2n+1 | : | 1 | 0 | 0 | 0 | 0 | 0 | 0 | 1 | 1 | 0 | 1 | 1 | 0 | 0 | 1 | 1 | 1 | 1 | 1 | 1 | 1 | 1 | 1 | 1 |
| | 3n+1 | = | 1 | 1 | 0 | 0 | 0 | 0 | 1 | 0 | 1 | 0 | 0 | 0 | 1 | 1 | 0 | 1 | 1 | 1 | 1 | 1 | 1 | 1 | 1 | 0 |
| | /2 | = | | 1 | 1 | 0 | 0 | 0 | 0 | 1 | 0 | 1 | 0 | 0 | 0 | 1 | 1 | 0 | 1 | 1 | 1 | 1 | 1 | 1 | 1 | 1 |

**Figure 1** shows 2 consecutive Collatz steps starting on an odd n0 shown in binary, an n0 that has several least significant 1 bits in a row before coming to a 0 bit. Each step ends when there is a divide by 2. The shortening of the string of least significant 1 bits is highlighted in yellow.

**Q2**: Is there anything predictable about how Collatz sequences fall into their final state of cycling between 1 and 4?

**A2**: This is another thing we readily see by applying the Collatz procedure in binary, i.e. that there are well defined exit points for any given sequence. I call these exits the ***picket fence*** numbers, because their appearance in binary resembles a picket fence, i.e. 1010101010….It is quite easy to see why these are the only ways to exit the Collatz sequence, using the same binary logic as above. Again, a multiply by 3 is equivalent to a multiply by two then adding the result to itself, i.e. 3n = 2n + n. So that, in binary, starting with a picket fence number, the Collatz procedure for the last odd integer in the sequence looks like:

| | | | | | | | | | | | | | | | | | | | | | | | | | | | |
|---|---|---|---|---|---|---|---|---|---|---|---|---|---|---|---|---|---|---|---|---|---|---|---|---|---|---|---|
| | n | : | | | 1 | 0 | 1 | 0 | 1 | 0 | 1 | 0 | 1 | 0 | 1 | 0 | 1 | 0 | 1 | 0 | 1 | 0 | 1 | 0 | 1 | 0 | 1 |
| + | 2n | : | | 1 | 0 | 1 | 0 | 1 | 0 | 1 | 0 | 1 | 0 | 1 | 0 | 1 | 0 | 1 | 0 | 1 | 0 | 1 | 0 | 1 | 0 | 1 | 0 |
| | 3n | = | | 1 | 1 | 1 | 1 | 1 | 1 | 1 | 1 | 1 | 1 | 1 | 1 | 1 | 1 | 1 | 1 | 1 | 1 | 1 | 1 | 1 | 1 | 1 | 1 |
| + | 1 | = | 1 | 0 | 0 | 0 | 0 | 0 | 0 | 0 | 0 | 0 | 0 | 0 | 0 | 0 | 0 | 0 | 0 | 0 | 0 | 0 | 0 | 0 | 0 | 0 | 0 |

**Figure 2** shows why picket fence numbers are always the last odd number in any Collatz sequence.

The result, of course, is a power of 2, which quickly reduces to 1 during the inevitable remaining steps. The first several picket fence numbers are shown in Table 1:

| Binary | Decimal | Decimal factorization | First # to exit | Number of n0 that exit* |
|---|---|---|---|---|
| **1** | 1 | 1 | 1 | 1 |
| **101** | 5 | 5 | 3 | 938,003 |
| **10101** | 21 | 3*7 | 21 | 1 |
| **1010101** | 85 | 5*17 | 75 | 23,743 |
| **101010101** | 341 | 11*31 | 151 | 37,687 |
| **10101010101** | 1365 | 3*5*7*13 | 1365 | 1 |
| **1010101010101** | 5461 | 43*127 | 5461 | 78 |
| **101010101010101** | 21845 | 5*17*257 | 14563 | 448 |
| **10101010101010101** | 87381 | $3^2$*7*19*73 | 87381 | 1 |
| **1010101010101010101** | 349525 | $5^2$*11*31*41 | 184111 | 36 |
| **101010101010101010101** | 1398101 | 23*89*683 | 932067 | 2 |

**Table 1**: The Picket Fence numbers, exit points for Collatz sequences. Picket Fence numbers divisible by 3 are excluded from Collatz sequences and hence are not connected to any other nodes. *Number of exits in the first 1,000,001 odd integers. The numbers sum to 1,000,001 – there are no other ways to exit the sequence.

Every third value in the picket fence sequence is divisible by 3, and hence not accessed as an exit point by Collatz sequences. In decimal notation, picket fence numbers alternate ending in 5 or 1 (it is easy to prove why this must be so), hence every other value is also divisible by 5. These facts will become important considerations later on.

It has been noted by many others that Collatz sequences, when looked at as graphs with every integer a node and places where different sequences intersect as branching points, tend to resemble trees when mapped together. We will use the graph/tree visual imagery extensively when discussing Collatz and Collatz-like sequences. For Collatz sequences themselves, if we exclude even numbers, the picket fence numbers described above can be thought of as the trunk of the tree, to which all of the major branches, sub-branches, and twigs attach themselves.

Exit points (the picket fence numbers on the trunk of the tree) are accessed less as they increase by some power law, and why this is so is, to a first order, fairly intuitive. In the first 50 positive, odd integers, n0 = 1 and 21 are only accessed as exit points by themselves, 21 because it is divisible by 3 and no other n0 can get to it. The next exit point (picket fence number) is at 85, and does not get accessed until n0 = 75. Thus all of the other odd numbers < 75 exit through the branch that attaches to the trunk of the tree at 5. And all of those numbers have nodes and branches attached to them, going up to infinity. Of all of the odd numbers ≥ 75, you might imagine that only a few of the remaining 13 that are < 100 attach to 85, and of the next 100 odd integers only a few % of those attach to the trunk at 85. A few percent will also attach to the next exit point at 341, since its first access reaches down to n0 = 151. Thus the exit points at 85 and 341 are attached to 2-4% of the nodes and branches as you go further up the number line. The next exit point is divisible by 3 and so cannot be reached, and then for the next exits (5461 and 21845), because they are so large and are already surrounded by a “thicket” of branches and nodes emerging from lower down the number line, only a few % of a few % of branches and nodes higher up will be connected to them, and so on. Obviously, there are second order effects, because the number of branches attached to 341 exceeds the number attached to 85, and this is contrary to the basic intuition that exit points lower down the tree should be connected to an ever-increasing number of branches further up, but to a first order the basic intuition holds.

**Q3**: Do Collatz sequences generally follow the ¾ factor reduction every 2 steps that is the default expectation, and that is why we think the sequences eventually always converge?

**A3**: We can test the average ¾ reduction hypothesis by building some null models:

1) The "low" sequence length model for $n0 = 2^k$, where the sequence descends straight to 1 in k steps because it consists only of 1/2 steps.
2) A "mid" sequence length model with ¾ drop every 2 steps, starting from an arbitrary n0, because we assume that 3/2 steps and 1/2 steps are equally likely on an average as the sequences descend.
3) A "high" sequence length model where you start with $n0*(3/2)^k$, where k is the bit length of n0. That is, you assume an initial run-up of k 3/2 increase steps as in $n0 = 2^k – 1$, then following the mid sequence length model after that, assuming a ¾ drop every 2 steps.

Results for some regular, exponentially increasing sets of n0 are shown in Figures 3 through 5. Generally, Collatz sequences do follow the mid sequence length, default model with 3/2 and 1/2 steps equally likely (Figures 4 and 5). As expected, we see some interesting behaviour in sequences starting with $n0 = 2^k – 1$ and, surprisingly, $n0 = 3^k$ (Figures 3 and 5). The amount of statistical persistence or autocorrelation (looks like clumpiness in the plots) of the $n0 = 3^k$ graph is almost disturbing (Figure 3c, Figure 5d, 5e, 5f). It is not a computational artifact, but is in fact a consequence of the branching process discussed in the previous section, as we shall see (Figure 6).

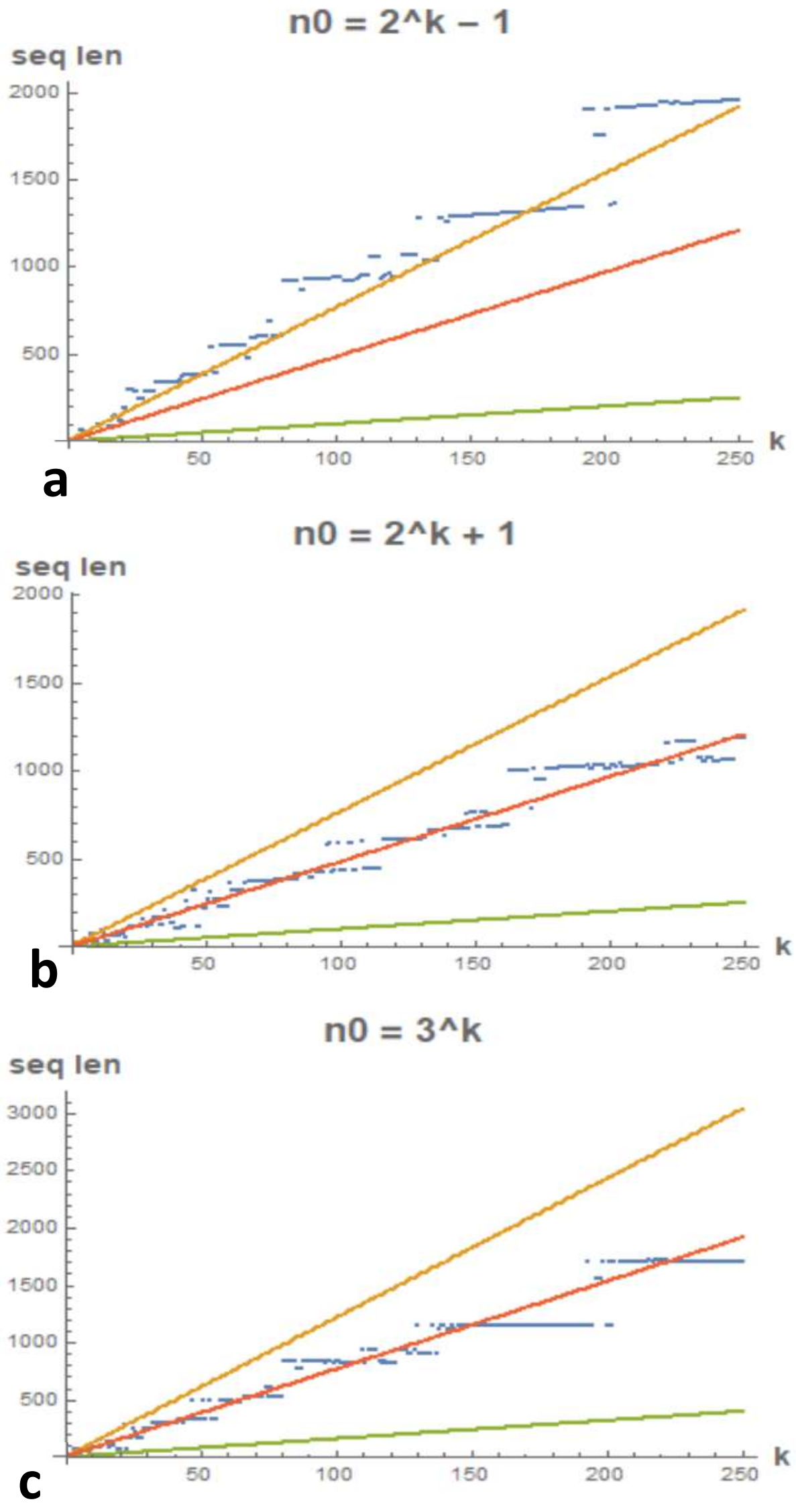

**Figure 3:** Plots for sets of exponentially increasing n0 according to the formulas specified in the title, where k is an integer increasing from 1 to 250. Actual sequence lengths are shown in blue dots. For n0 = $2^k - 1$ (subplot **a**), the high null model is followed (amber line) as expected; for n0 = $2^k + 1$ and n0 = $3^k$ (subplots **b** and **c**), the mid sequence length prediction is followed (red line), also as expected. What is not expected is the high degree of persistence in the plots. The low sequence length model is shown in green.

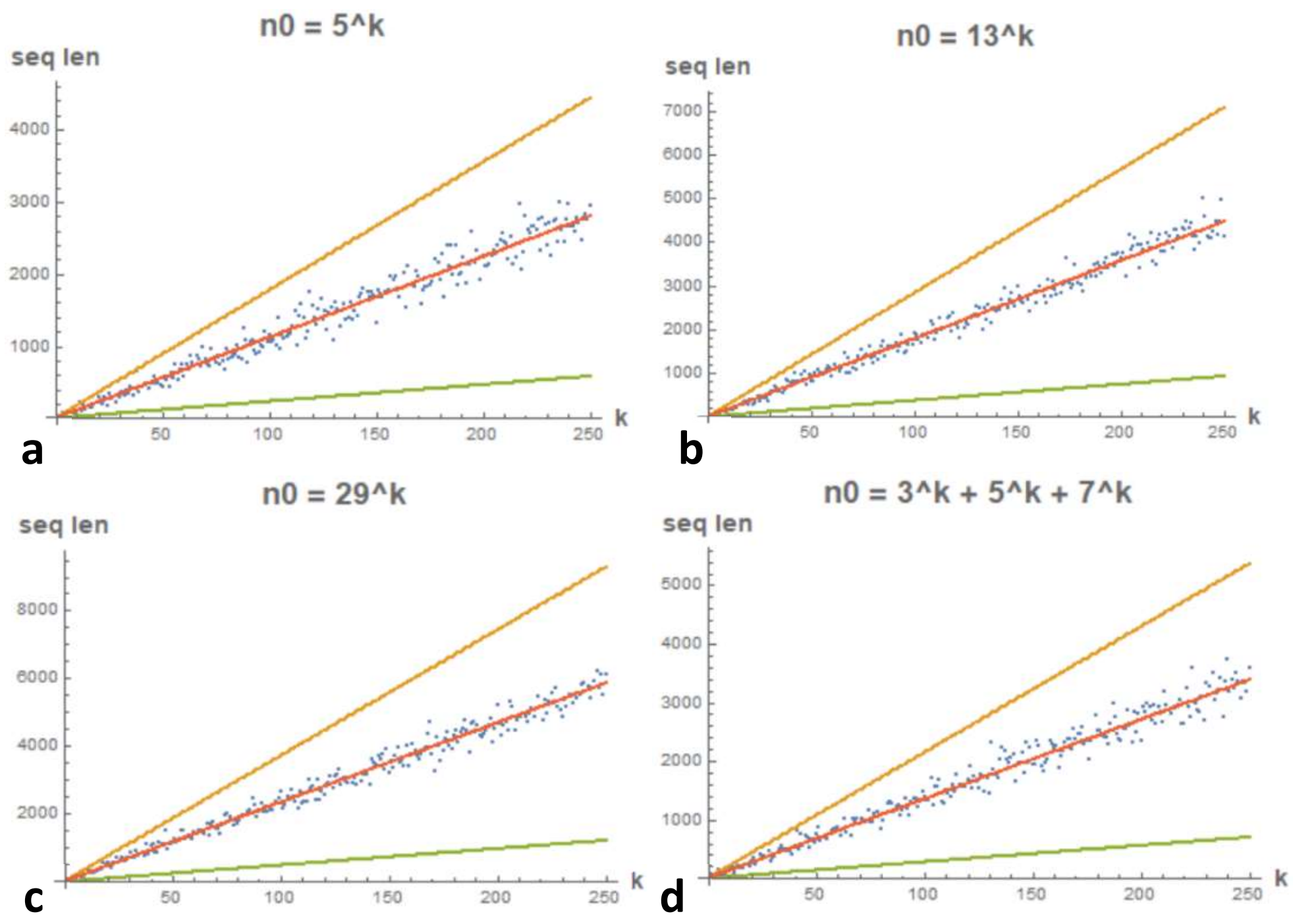

**Figure 4:** Plots for various sets of n0 = $p^k$, where p is a prime number raised to an integer power k. The color scheme for the plots is the same as Figure 3, e.g. actual sequence lengths are shown with blue dots, and the high, mid, and low sequence length models are shown with amber, red, and green lines.

We can see that n0 in the series $2^k - 1$, as expected, either meet or exceed the high model (mostly exceed, Figure 3a, Figure 5a). Like the L.A. weather, they are also quite persistent. The series $2^k + 1$ sticks to the mid sequence length null model, as does $3^k$, but both have large "islands of persistence" or near-persistence for reasons that are not immediately obvious (Figure 3b,c and Figure 5b,d). After that, similar series of n0 = $p^k$, or n0 = $p^k + p^j$ +..., but with p = prime numbers > 3, stick to the mid sequence length null model quite well overall (Figure 4).

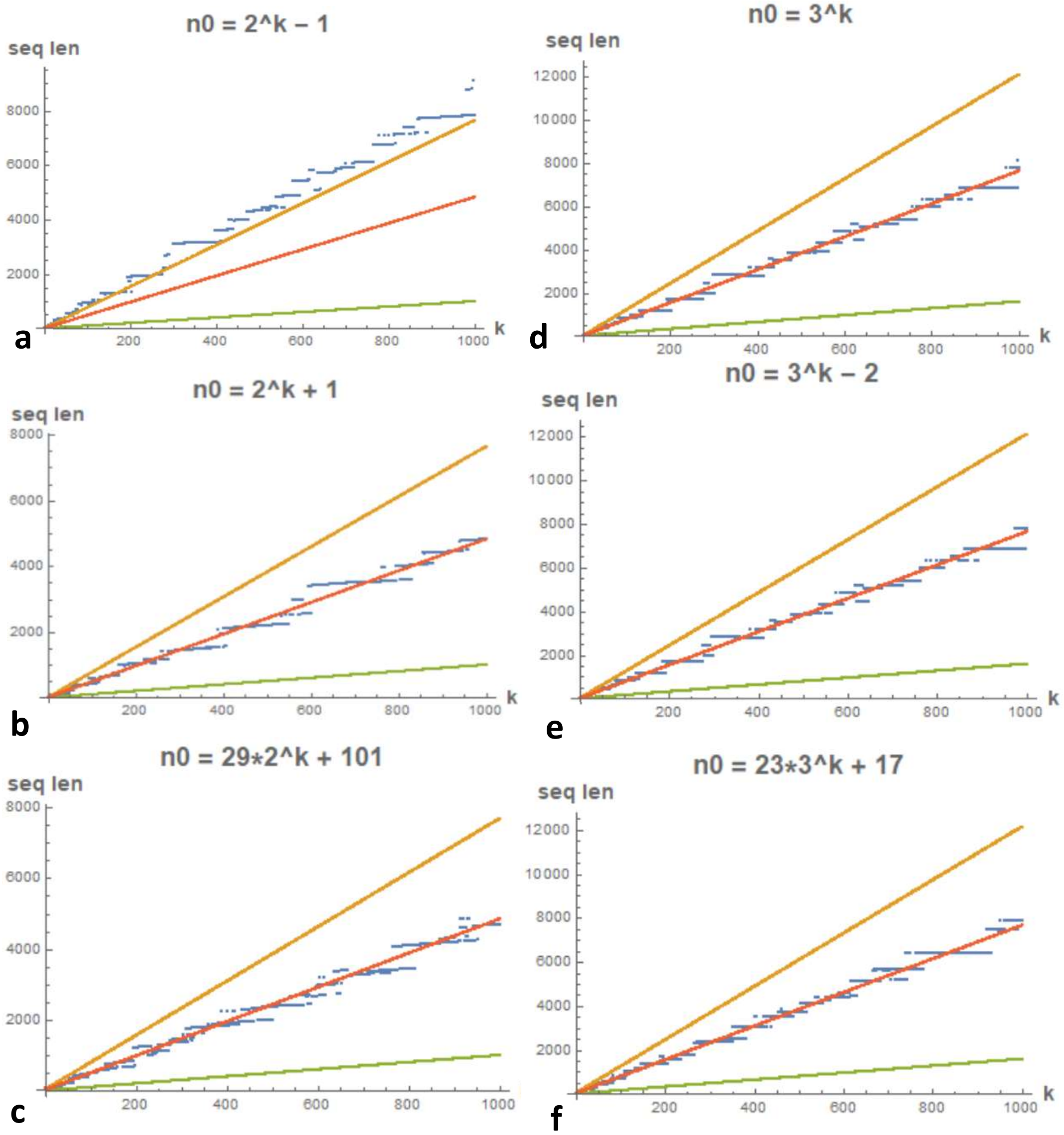

**Figure 5:** Plots for Collatz sequence length for sets of n0 similar to Figure 3, but for much larger k. The persistence of all of these plots persists, even when modified by various offsets and multiples!

We can explore the persistence/near-persistence of the n0 = $2^k$ and $3^k$ sequences by extending them further, to k = 1000 (Figure 5). Both show large islands of persistence in sequence length, but the largest are found in the series for n0 = $3^k$ and related. For n0 = $3^k$, between k = 869 – 981 inclusive, the sequence length (using our definition of Collatz sequences here) is exactly 6842, except for seven values from k = 886-892 where sequence length = 6311, and one value at k = 971 where it is 7804. This even though n0 increases by 53 orders of magnitude over that interval! It's tempting to think that the plot for n0 = $3^k$ does this as a consequence of the fact that powers of 3 are excluded from the Collatz sequence generator. But the plot for n0 = $3^k$ – 2 (subtracting 2 keeps the n0 an odd number, Figure 6d) looks almost identical, even though none of the n0 = $3^k$ – 2 are divisible by even a single factor of 3. Going from n0 = $2^k$ - 1 to $2^k$ + 1 has a very large effect on sequence length although, oddly enough, for n0 = $2^k$ + 1 the results are even more statistically persistent than they are for $2^k$ - 1.

We can explore these islands of persistence and how they happen in Figure 6. Figure 6 examines a large island from Figure 3c, between the values of k = 130 - 200. The log base 10 of each successive number in the sequence is plotted vs the i'th iteration, until it reaches n = 1 where the sequence terminates. We can see in blue, amber, and green three sequences at the beginning, middle, and end of the island, all converging to the same sequence branch of the Collatz tree. That they all have the same sequence length of 1144 iterations is a strong clue that they are on the same branch. For example, another branch, in red, appears on the plot to be on the same island but in fact converges in 1155 iterations. Still, the statistical persistence of n0 = $3^k$ is a profound mystery. We haven't proved that islands of persistence are 100% composed of sequences on the same branch, but the fact that they have the exact same sequence length, together with plots like Figure 6 verifying that they are on the same branch, is highly suggestive.

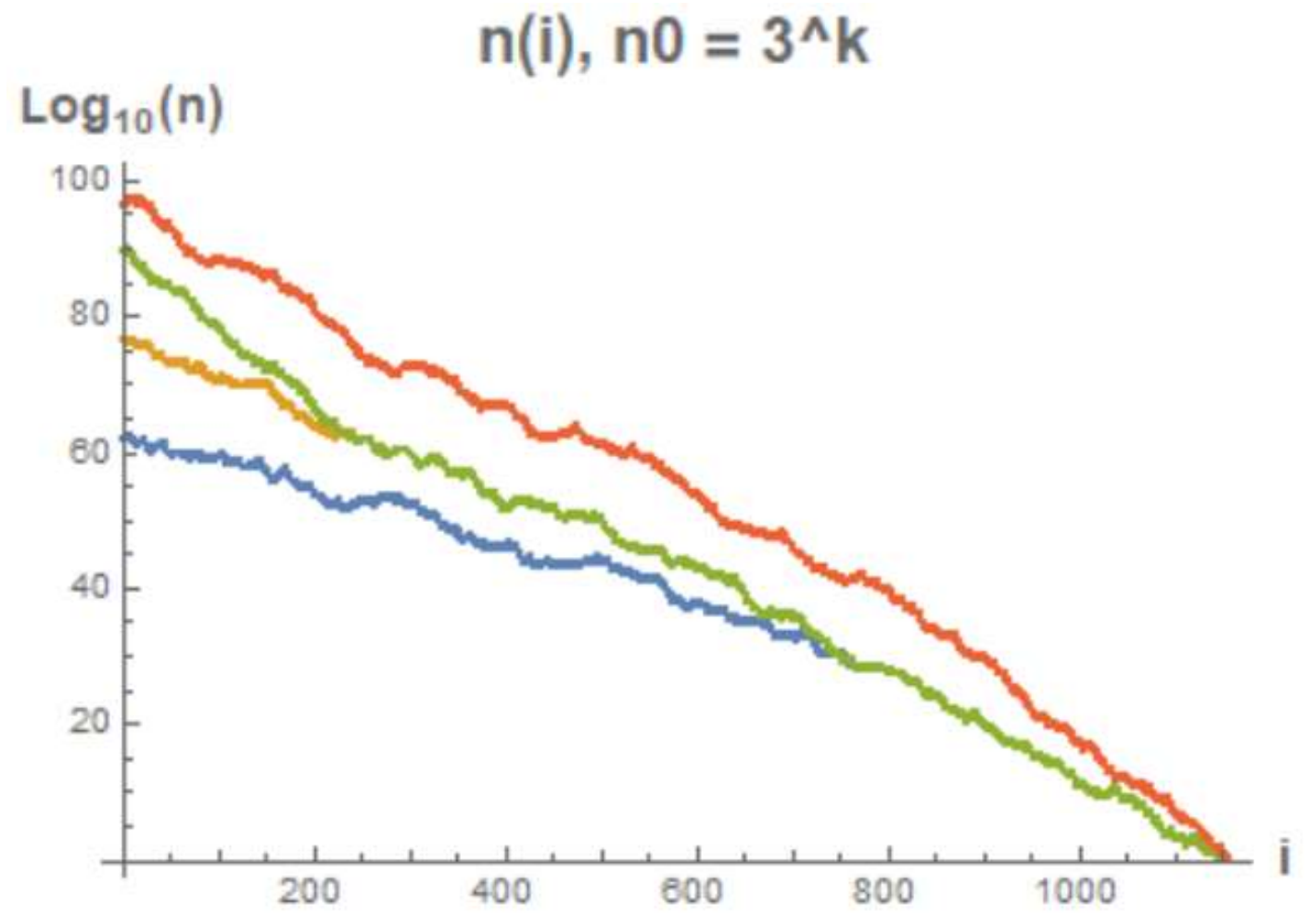

**Figure 6:** An actual set of sequences from an island of persistence in the set n0 = $3^k$. n = n0 at i = 1, until the sequence converges to 1 in the Collatz program. In blue, amber, and green are the sequences for k = 130, 160, and 188 at the beginning, middle, and end of the island of persistence, all converging in exactly 1144 iterations. In red is the sequence for k = 202, which looks in Figure 3c to be on the island (or the archipelago), but converges in 1155 iterations instead.

For our purposes, the large islands of persistence serve as a type of x-ray to reveal the many hundreds of sub-branches connected to any given major branch of the Collatz tree. Horizontal cross-sections of the tree structure at constant sequence length reveal that the branches are spaced logarithmically in powers of 2 and 3. Given that the islands still form with arbitrary shifts and multiples of the set of exponentially spaced n0 (Figure 5c and 5f), it shows that all n0 are most likely well connected to many others in these branched structures. This makes the idea that a cycle or loop of integers > 250,000 long and sitting high up the positive number line, as per previously determined limitations on the Collatz conjecture, is not connected to any of the thicket of branches rising up from below increasingly unlikely. I imagine that, if islands of persistence for sets of sequences beginning n0 = $3^k$ and n0 = $2^k$ + m could be proved to continue indefinitely, proof of the Collatz conjecture itself would be a straight forward consequence.

**Q4**: What makes Collatz sequences interesting and how would we define or build other types of similar sequences that are interesting?

**A4**: First we will define types of sequences as interesting if:

1) They don't appear to ever diverge. Note that divergence at only a finite number of n0 is not really a coherent possibility given the way Collatz sequences branch upwards and how they work in general. If one n0 diverges then an infinitely many will, although they could co-exist with infinitely many that do converge, as we shall see below.
2) They rise up and down a lot, i.e. many of the sequences don't tend to converge monotonically.
3) There is a large variation in sequence lengths, which are not correlated -- except to a first order -- to the magnitude of the starting points, n0.
4) They don't always converge to the same point, but have several stable points of convergence, corresponding to "loops" of integers that vary widely in length and height above one. It is not known if the number of stable points can be infinite, but usually they appear to be finite and quite small in number (< 10).

Collatz sequences for a variety of n0 are interesting in that they appear to have the first three qualities but not the fourth: they only appear to have one very short loop that bottoms out at n = 1, as is stated in the original Collatz conjecture. If we include the simple modification to the Collatz of 3n – 1 (subtract one instead of add), there are three loops that bottom out at 1, 5, and 17 (Lagarias 1985). Below we see that there are Collatz-like sequences that exceed Collatz sequences in the first three qualities while also having the fourth. Some of these are closely related to Collatz sequences, and demonstrate again (i.e. as other authors have done in different ways) why the Collatz conjecture is very likely to be true.

### Alternatives to Collatz sequences and their characteristics

By adding an ElseIf statement to the original formulation of the Collatz conjecture (which has only an "If-Then-Else" structure), we can create alternate sequence generating programs, some of which are closely related to the original Collatz sequence generator. We can use simple probability calculations to build null models for how we expect these sequences to fall towards lower integers from a given n0, and how we expect the different trees that are produced to branch, and see where those expectations are met or contradicted. As mentioned above, the alternate programs do not always converge to 1, and show some interesting behaviours which exceed the original Collatz formulation in doing "interesting" things, but also suggest why the Collatz conjecture is likely true. The alternative programs are labeled P#, where the "#" indicates some integer (# increasing in the order of their discovery). Only the more interesting and Collatz conjecture relevant programs are discussed below. One thing we shall see is that there appears to be only a finitely many formulations that remain bounded and don't diverge in any of their sequences as far as we can detect. Also, where the alternative formulations differ from the original Collatz, the parameters that define the differences are not large, i.e. there don't appear to be any cycles or loops with millions of iterations that start at large numbers many dozens of digits long.

**P2**: P2 is a Collatz-like set of sequences that tinker with the original rules in ways that it was hoped would create more meta-stable sequences, i.e. sequences that hopped around more wildly and unpredictably than Collatz sequences, yet remained stable and converged eventually. That is, it was hoped that slower average convergence than the expected ¾ factor every 2 steps in the Collatz procedure would make for even greater variation in sequence length, and find some very interesting behaviours according to all four criteria for what is "interesting". P2 is defined like so:

If even: /2; ElseIf divisible by 3: (7n/3 + 1)/2; Else: (5n + 1)/2

P2 is different from the original Collatz conjecture (which we will call P1) in 2 ways: it distinguishes between odd integers divisible by 3 and all other odd integers, and it uses different multiples for the odd integers. P2 produces eight known cycles or complete loops of integers, as shown in Table 2 below.

| Lowest n | factors | # iterations | Lowest root-node | factors | Highest n | % of exits |
|---|---|---|---|---|---|---|
| **1** | 1 | 4 | 1 | 1 | 4 | 10.94 |
| **7** | 7 | 6 | 7 | 7 | 28 | 16.04 |
| **21** | 3*7 | 310 | 5 | 5 | 16443858 | 58.06 |
| **85** | 5*17 | 6 | 85 | 5*17 | 340 | 1.96 |
| **121** | 11*11 | 6 | 113 | 113 | 354 | 10.21 |
| **141** | 3*47 | 6 | 77 | 7*11 | 564 | 1.37 |
| **1303** | 1303 | 33 | 521 | 521 | 53764 | 1.31 |
| **69721** | 113*617 | 44 | 20981 | 20981 | 4228008 | 0.11 |

**Table 2**: A comprehensive list of all known loops in P2 and some of their characteristics.

Table 2 also shows some of the characteristics of the loops in P2, and reads like so, from left to right: the lowest value in the repeating set of integers, the prime factors of that value, the number of iterations the program goes through before it repeats itself, the lowest integer that gets caught in the loop when it starts iterating in the program, the prime factors for that lowest integer, the highest integer in the loop, and the percentage of time that a random starting number will wind up falling into that loop. Although the "% of exits" are assessed for the first 500,000 odd integers, tests on random integers show that they exit into the P2 loops in the percentages shown.

There are some interesting things to note about the loops produced by P2 as shown in Table 2. The most obvious is the large loop starting at n = 21, with 310 iterations that rise quite high before eventually coming back to 21, its lowest point. The loop's (hence forth "L21" for P2) boundaries contain all of the higher loops in Table 2, and it has "roots" that reach down to the integer 5 that feeds into it. Because the loop starts low and has so many iterations, it is connected to the highest percentage of nodes and branches above it, as shown in the right most column with 58%. Curiously, there was no prominent exit at n = 21 in the original Collatz sequence generator, as the numbers containing a factor of 3 cannot be embedded in any sequence.

**Null model for P2 sequence lengths**

If we compute a very simple null model for P2, similar to the factor of ¾ drop every 2 steps for the original Collatz, we can calculate that we would see a decrease by a factor of 0.91 in the sequence every 6 steps – quite a bit slower than the original Collatz – like so:

$0.91 = (1/2)^3(7/6)(5/2)^2$

Thus, we predict that the sequences emerging from P2 are only barely convergent. This makes some intuitive sense, for most of the time, on encountering odd numbers, the sequence is getting a 5/2 multiple, not a 3/2 as in the original Collatz program. If it were not for the intermediary step for odd integers divisible by 3, the sequences would tend to diverge. If we plot out actual P2 sequences against this simplistic null model, however, we see that the P2 sequences drop quite a bit faster than that simple null model. But there is an easy, intuitive fix to the simplistic null model which improves the model performance a lot: by recognizing that P2 tends to exclude integers divisible by 5 due to the 5n + 1 step (similar to the way the Collatz program excludes factors of 3), we see that factors of 3 in the odd integers of the sequence are somewhat enriched due to the near absence of the factors of 5. So if we enrich the factors of 3 in the sequence by 25% (i.e. multiply by 5/4 due to the near absence of numbers that are multiples of 5), we get an 0.62 decrease in 10 steps.

$0.62 = (1/2)^5(7/6)^{2.08}(5/2)^{2.92}$

This might not seem like much of a difference at first glance, but the new model results in predicted sequence lengths that are almost a factor of 3 shorter than the simplistic null model. Still, this null model descends considerably more slowly than the original Collatz procedure, which declines by an average factor of 0.24 in 10 steps. This modified null model performs surprisingly well, although it still drops slower than the actual sequences even with the assumption that all integers containing a factor of 5 are excluded (Figure 7). In fact, multiples of 5 are not excluded from the sequence by definition, they just appear much more infrequently, because the 7/3 step allows the program to produce a factor of 5 sometimes. Still, if we check sequences from P2 we see that the modification to the null model is mostly justified: odd integers containing a factor of 3 are enriched compared to the 33% we would expect at random, and factors of 5 are much reduced. We see almost exactly the 25% enrichment for multiples of 3 we had calculated, to about 40% of all odd numbers encountered (in randomly selected samples of sequence), and only 5% of the odd numbers are factors of 5. This is what allows the P2 sequences to converge as quickly as they do. This easy modification to the simplistic null model doesn't always work for other programs, however, as we shall see.

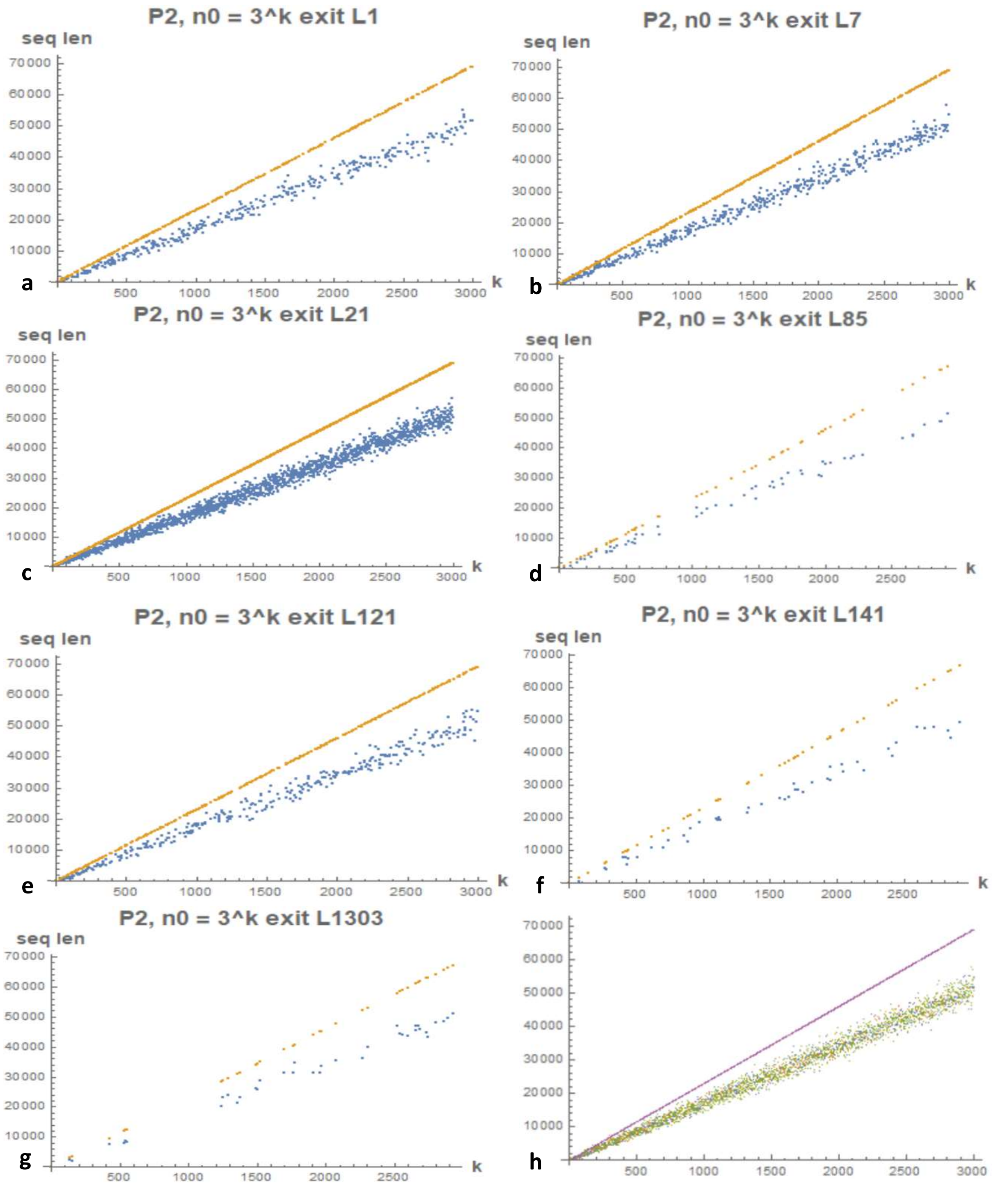

**Figure 7:** Comparison of actual sequence lengths for P2 with sequence lengths predicted by the modified null model, starting from the set of integers n0 = $3^k$, for k = 1-3000. Sequence lengths shown are from n0 to the lowest point in the loop in which they exit (blue dots, sub-figures a – g). The modified null model of an average factor of 0.62 decline every 10 steps of sequence is shown in gold (a – g). The plot for L69721 is not shown as it contains only 2 points (roughly in keeping with the expectation from Table 2). Points exiting all loops are overlaid in sub-figure h, with the null model from L21 in purple.

**Null models for exit points and branching**.

One can gain some intuition as to why, for example, L1 in P2 reaches so many more branches and nodes than L69721, as shown in the leftmost column of Table 2, with a simplistic calculation. Assuming L1 connects to 11% of all nodes above it on the way up to L69721, it will be connected to about 7700 nodes as you go through the first 70,000 n0. At about 70,000 on the positive number line, L69721 is only connected to the 44 nodes in its loop. That ratio, 44/7700, is not too far from the actual estimated ratio of nodes that each loop is connected to as you travel further up the number line. Of course, this is an inexact calculation and meant for rough illustration for how frequently each loop connects to the numbers (nodes) above it.

This idea that the nodes in any given loop are equally likely to connect to the nodes and branches above them needs some modification, however. The loops at 85, 121, and 141 all consist of shuffles of the same 6 operations: two divide by 2 steps, three divide by 3 steps, and one multiply by 5 step (note that 1/2*1/2*7/6*7/6*7/6*5/2 = 0.9925, which helps explain why this simple combination is so rich in loops). Yet L121 leads to many more branches above it than either of the other 2 loops, for reasons that are not immediately obvious. The cycles of integers, shown below, suggest some reasons:

85 -> 213 -> 249 -> 291 -> 340 -> 170 -> 85
121 -> 303 -> 354 -> 177 -> 207 -> 242 -> 121
141 -> 165 -> 193 -> 483 -> 564 -> 282 -> 141

L85 contains three numbers that are multiples of 5, and therefore less accessible to sequences descending from above, as we discussed in the section on null models for P2. Only two of the numbers in L85 are divisible by 3, whereas four of the L121 sequence are, and these are also more accessible. L141 has many of its nodes higher than the rest, and so tends to miss the critical nodes lower down the number line. Both L85 and L141 have consecutive ½ steps, and these lead to fewer node and branch connections.

**P4**: P4 represents a set of closely related, Collatz-like programs that intersect the original Collatz sequences when the parameter m=15, and demonstrate through their stability well past m=15 that the Collatz sequences are probably always convergent. These programs are equivalent to the original Collatz except when they encounter an integer divisible by 5, at which point they divide by 5 and multiply by an integer (labeled “m” throughout the rest of this paper), before adding 1 and then dividing by 2 as usual. Because 15/5 = 3, when m = 15 the program is equivalent to the original Collatz program. Table 3 shows a compendium of interesting loops and other characteristics for these simple programs at various values for m.

While most m in Table 3 are primes, some m that are not prime were worked up for the number and location of loops for various reasons. The progression of programs with m = 7, 15, 31 all having only one loop was suggestive of a pattern behind why there might be only a single loop associated with some programs, and with the original Collatz itself. Alas, m = 63 broke the pattern, as it has 4 loops.

**P4:** If even: /2; ElseIf divisible by 5: (m*n/5 + 1)/2; Else: (3n + 1)/2;

| m | Loop lowest n | # iterations | Loop highest n, >20 iter | Note |
|---|---|---|---|---|
| **7** | 1 | 2 | | |
| **11** | 1, 3 | 2, 3 | | |
| **13** | 1, 5, 23, 25 | 2, 8, 3, 3 | | Only example of loops for adjacent odd numbers 23, 25 not including 1 |
| **15** | 1 | 2 | | Equivalent to the original Collatz |
| **17** | 1, 5 | 2, 10 | | |
| **19** | 1, 5 | 2, 2 | | |
| **23** | 1, 3, 23, 65 | 2, 4, 7, 6 | | |
| **29** | 1, 5, 445, 595 | 2, 10, 13, 13 | | |
| **31** | 1 | 2 | | |
| **35** | 1, 5, 23, 43, 365, 9355 | 2, 11, 12, 6, 25, 156 | n0 = 365: 15,572<br>n0 = 9355: 77,849,600 | |
| **53** | 1, 25, 35, 43, 55, 63, 2125, 15871 | 2, 5, 5, 199, 5, 5, 814, 179 | n0 = 43: 4,239,444<br>n0 = 2125: 946,605,753,764,320<br>n0 = 15871: 45,323,252 | Longest loop discovered at n0 = 2125, has 814 iterations |
| **55** | 1, 5, 263 | 2, 11, 49 | n0 = 263: 3,402,836 | Appears to have an odd # of loops |
| **59** | 1, 5 | 2, 20 | | |
| **61** | 1, 5 | 2, 139 | n0 = 5: 5,561,216 | |
| **63** | 1, 95, 191, 203 | 2, 10, 20, 10 | | |
| **67** | 1, 5, 925 | 2, 7, 224 | n0 = 925: 15,482,864,544 | Also an odd number of loops |
| **73** | 1, 5, 71, 505 | 2, 12, 59, 808 | n0 = 71: 129,032<br>n0 = 505: 430,352,418,705,748 | Some sequence lengths 10's of millions digits rising to >10^1300, yet starting from n0 < 1000 |
| **101** | 1, 5, potentially others | 2, 14, others | | Most n0 may diverge, but there are infinitely many that converge |
| **127** | 1, | | | Most n0 likely diverge, but infinitely many known to converge |

**Table 3** A compendium of all known loops in P4, for different values of m, with the highest integer in the loop listed only for loops > 20 iterations. Values of m not on the list were never tried by the author.

### Null models for P4 sequence lengths

A simplistic null model for the average rate of decrease in the sequences of P4, similar in formulation to P2, would be:

$\eta_1 = (1/2)^5(m/10)^1(3/2)^4$

That is, we expect to see a factor of $\eta_1$ decrease in 10 random steps because half of the integers encountered will be even, on average, and of the odd integers we expect one in five to be divisible by 5.

If we try to modify the null model a bit like we did with P2, i.e. accounting for some enrichment in factors of 5 in the sequence (hence more m/5 steps) and a weeding out of factors of 3 due to the (3n + 1)/2 step, we would get:

$\eta_2 = (1/2)^5(m/10)^{1.35}(3/2)^{3.65}$

If we set $\eta_1 = 1$ to determine the boundary of stability, i.e. the value of m when the sequences start to become unstable and diverge, we get m = 63.21. On the other hand, setting $\eta_2 = 1$ results in m = 43.53. As we see in Table 3 and we shall see more below, in terms of this boundary of stability and how the sequences actually descend, this time the very naïve null model $\eta_1$ is closer to the truth.

We explore the $\eta_1$ null model in Figures 8 and 9, below, for the particular case of P4 m=53. This program has eight exit loops, like P2, listed in more detail in Table 4 below. Figure 8 explores sequence length for the first 10,000 positive, odd integers, hence the x-axis is linearly spaced. Figure 9 explores sequence length for the first 1,000 powers of 3, similar to Figures 3c and 5d,e,f, hence the ticks are exponentially spaced. Although it appears in Figure 9 that, like P2, our null model underestimates the average rate of decrease for the sequences in P4 m=53, we see in Figure 8 that is not quite correct. L43, L55, L63, and L2125 all have significant sets of sequences that far exceed the null model, which is not apparent in Figure 9. Why this is so is not clear. None-the-less, P4 m=53 faithfully converges to one of the eight exit loops in all of the many 1,000's of cases across 100's of orders of magnitude that have been tested so far.

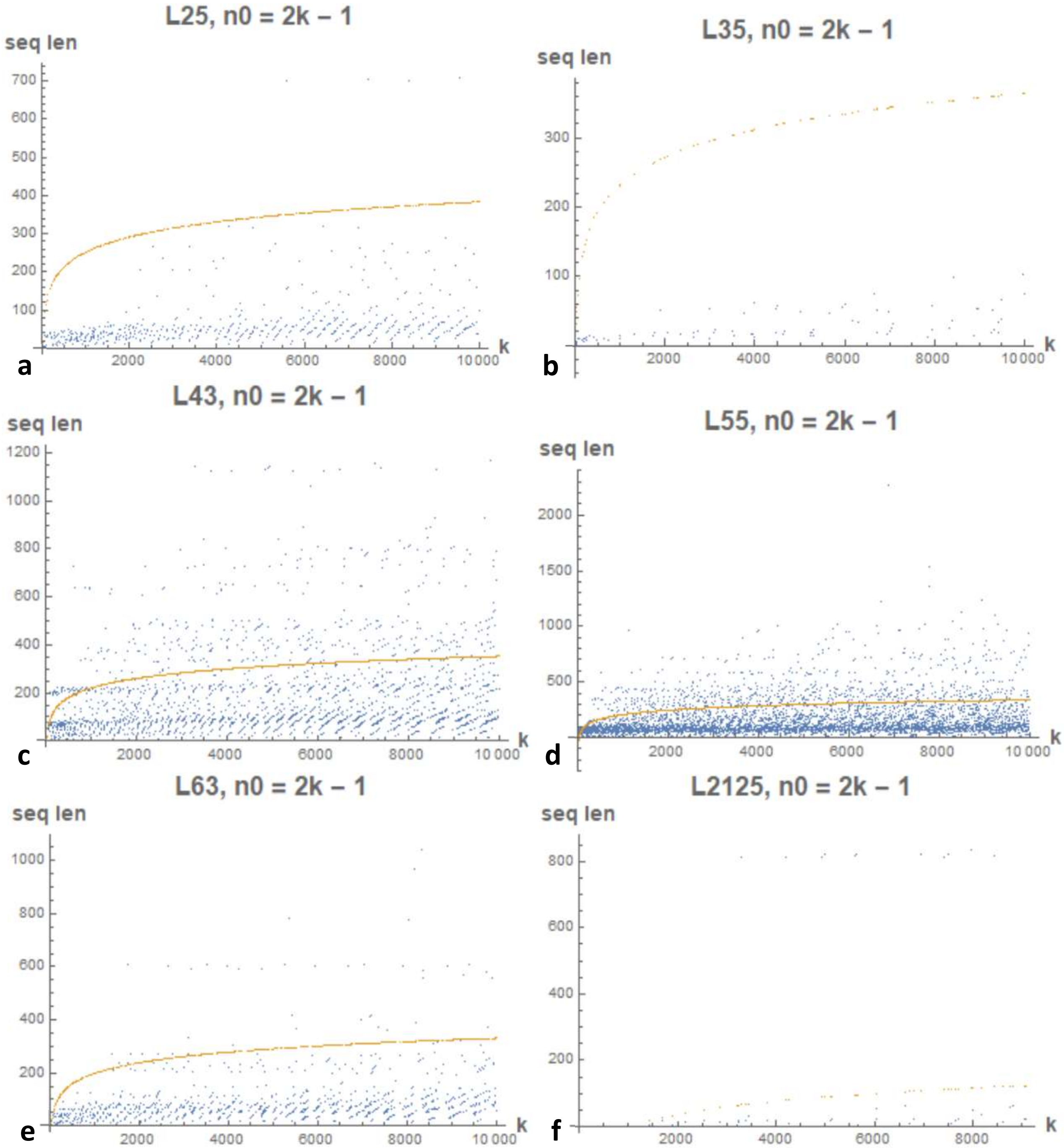

**Figure 8:** Plots for the sequence lengths of the program P4 m=53 for n0 = {the first 10,000 odd integers}. Only 6 of the 8 exit loops in P4 m=53 are shown for clarity: the subplot for L1 looks similar to L35, with fewer sequences and the sequence lengths quite low, and the subplot for L15871 is very sparse with only a handful of points, as expected. Sequence lengths are indicated by blue dots, the corresponding simple null model prediction for each point is shown in gold. Note the large variation in y-axis scale.

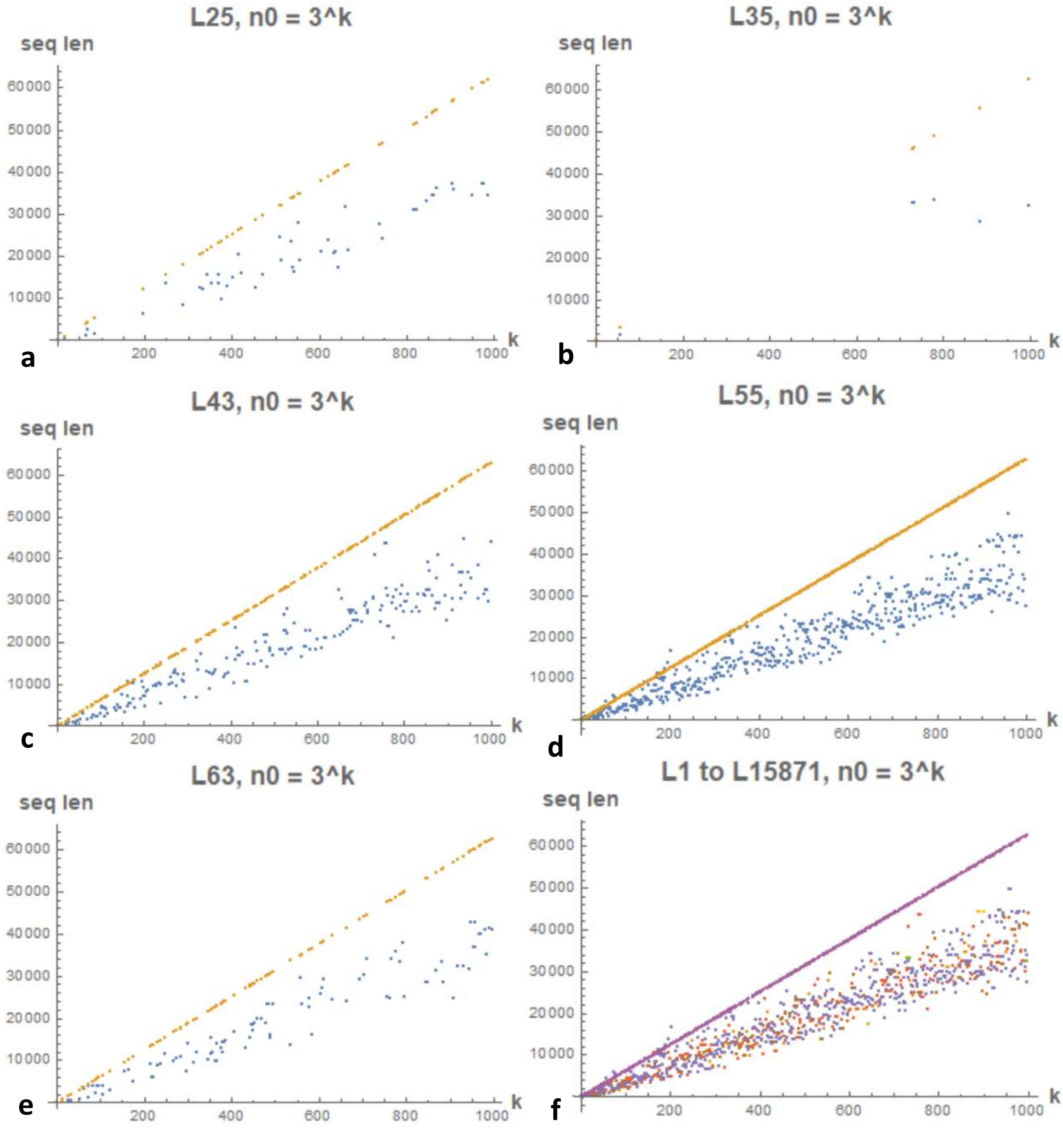

**Figure 9:** Similar to Figure 8 for P4 m=53, but for a different set of n0, where n0 increase exponentially as n0 = $3^k$ for k = 1…1000. Actual sequence lengths are marked with blue dots; the corresponding simple null model predicted lengths in gold. The subplot for L1 had no points in this set of n0 and was excluded; subplots for L2125 and L15871 were very sparse, resembled the subplot for L35, and were excluded for brevity. Points exiting all 8 loops together are shown in subplot f, with the null model for L55 shown in purple.

**Null models for branching and exit points for P4**

**P4, m=53** If even: /2; ElseIf divisible by 5: (53n/5 + 1)/2; Else: (3n + 1)/2

| Lowest n | factors | # iterations | Lowest root-node | factors | Highest n | % of exits |
|---|---|---|---|---|---|---|
| **1** | 1 | 2 | 1 | 1 | 2 | 0.16 |
| **25** | $5^2$ | 5 | 25 | $5^2$ | 200 | 5.73 |
| **35** | 5*7 | 5 | 23 | 23 | 186 | 0.62 |
| **43** | 43 | 199 | 43 | 43 | 4,239,444 | 23.90 |
| **55** | 5*11 | 5 | 3 | 3 | 292 | 58.4 |
| **63** | $3^2 7$ | 5 | 63 | $3^2 7$ | 504 | 10.0 |
| **2125** | $5^3 17$ | 814 | 2125 | $5^3 17$ | 946,605,753,764,320 | 0.43 |
| **15871** | 59*269 | 179 | 8359 | 13*643 | 45,323,252 | 0.15 |

**Table 4** A list of all of the known loops in the P4 program with m= 53, similar to Table 2 for P2. Exit % for the first 50,000 odd integers > 20,000 (i.e. above the last loop).

Table 4 has some remarkable features. The loop at 2125 rises higher and iterates longer than any other loop in any of the programs discovered so far. The loops L25, L35, L55, and L63 all contain permutations of the same 5 operations and cover roughly the same small domain on the positive number line, and yet they behave so differently.

25 -> 133 -> 200 -> 100 -> 50 -> 25
35 -> 186 -> 93 -> 140 -> 70 -> 35
55 -> 292 -> 146 -> 73 -> 110 -> 55
63 -> 95 -> 504 -> 252 -> 126 -> 63

In particular, L55 has 2 ***orders of magnitude*** more branches and nodes connected to it than L35, even though superficially both loops of integers are very similar. The loop L43 has 199 iterations and, by the logic of the analysis of P2, should have many more branches and nodes than L55, which only has 5 iterations, but in fact has only ~40% as many. Actually, the eight loops of P4 m=53 don't really behave like the eight loops of P2 at all when it comes to predicting how many branches and nodes they will be connected to. P2 is much more intuitive, and tends to follow the logic that any given node in any given loop is on average equally likely to lead to branches and nodes further up the tree.

We can explore a little what is going on with L35 and L55 by plotting the sequence trees of the first 10 positive integers that fall into those loops, as shown in Figure 10. Some trends are immediately apparent: the L35 sequences are much shorter than L55 (as can also be seen in Figure 8 b,d). Shorter sequences mean fewer attachment points for other sequences descending from up above. Also, L55 has deep roots: of the first ten positive, odd integers, nine of them fall into L55. The first ten integers to exit at L55 don't even make it half way up the number line to the bottom of L55, whereas eight of the first ten integers to exit L35 are greater than 35 – some much greater. As we saw with the original Collatz program in Table 1, more nodes on the tree lower down the number line means many, many more nodes to connect with further up.

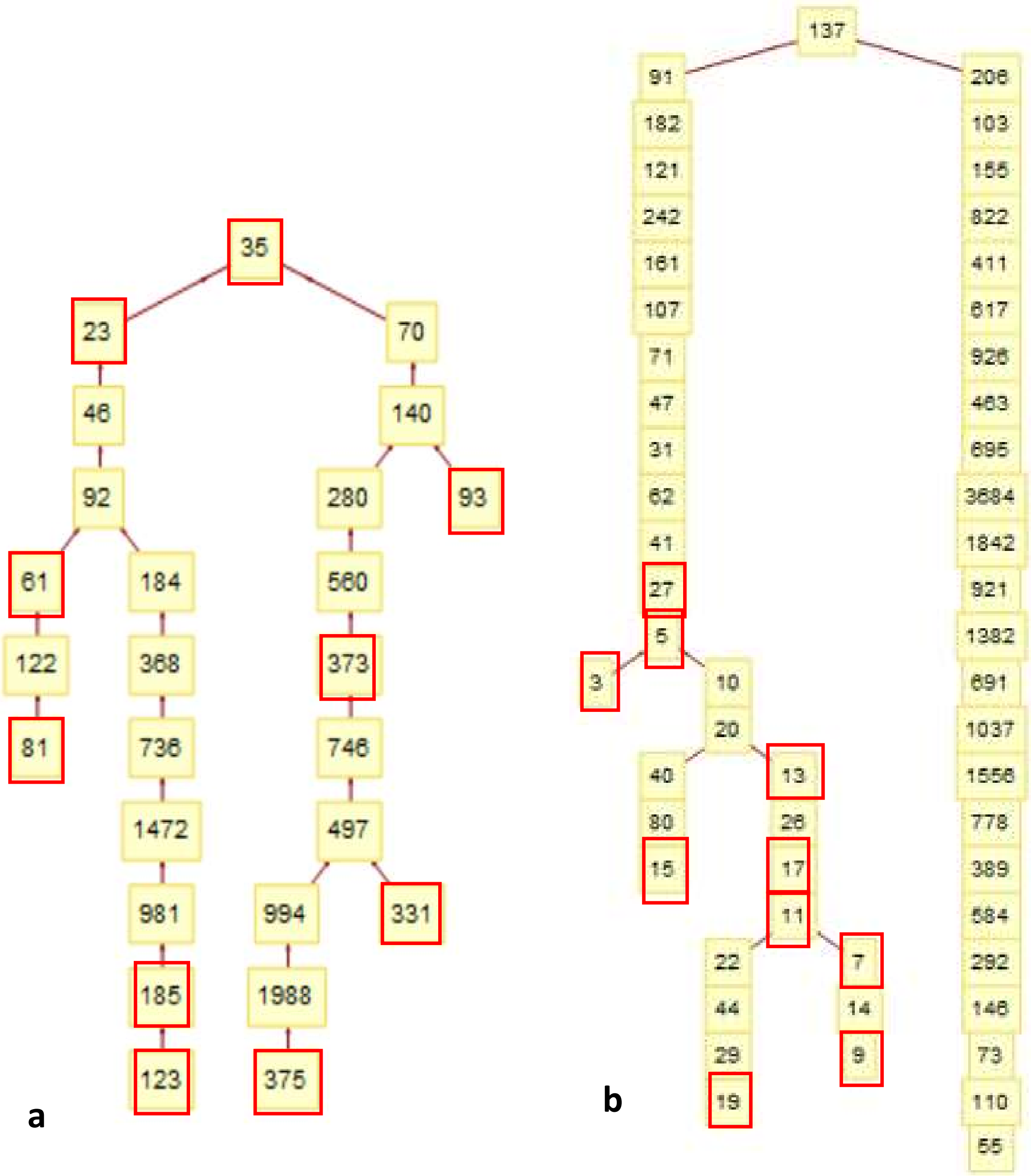

**Figure 10:** Tree plots of the first 10 odd integers to exit L35 (a) or L55 (b) for P4 m=53. The first 10 n0 for each tree are outlined in bright red.

**P4 m=73**

P4 is equal to the original Collatz sequence generator at m = 15, so as m increases from that point one wonders what happens to the stability of the sequences, at what point they start to diverge. In fact, as we saw in Table 3, m can increase to nearly 5 times more than m=15 and still remain convergent for all sequences as far as we can see. At m = 73, however, the sequences for some n0 appear to have metastable properties, going on for millions of iterations and rising to huge numbers, even though they eventually converge. Table 5 shows some examples of these n0 that have been discovered so far for P4 m=73. It is not known whether the last n0 entry ever converges or not. It is the only n0 in many thousands of random trials, some which rise to millions of iterations, for which I cannot verify convergence and, in fact, seems to diverge. The odd numbers to each side of the last n0, i.e. n0 = 20,387…+/- 2, both converge quite well (in 5388 and 1677 iterations and exiting L5 and L71, respectively). There is something very unusual about that exact n0. Whether P4 m=73 converges eventually for all n0 is an open question. If it does not, the next possibility for the boundary of stability in P4 is m = 71, which is somewhat suspicious.

| n0 | # of iterations | Exit loop | Max # in sequence (approx.) |
|---|---|---|---|
| **665** | 7,052,259 | L5 | $\sim 10^{680}$ |
| **26301 and 26303** | 7,123,078 | L5 | $\sim 10^{680}$ |
| **24,256,988,407,731,586,453,766** | 8,994,784 | L71 | $\sim 10^{1300}$ |
| **24,302,489,955,379,179,567,675** | 16,247,650 | L71 | $\sim 10^{1600}$ |
| **8,289,328,408,488,871,421,602** | ~19,200,000 | L71 | $\sim 10^{1300}$ |
| **20,387,503,213,950,575,404,289** | >140,000,000 | ? | $>10^{8336}$ |

**Table 5** Some examples of the longest discovered sequences in P4 m=73.

**P6**: A Collatz-like closely related set of programs that each exclude a prime integer from their sequences, similar to the way the original Collatz excludes factors of 3, but are engineered to be convergent. These Collatz-like sequences are built by adding extra ElseIf-statements to the original program to try and slow down the explosive growth in the sequence if the original 3/2 steps were replaced by p/2 steps, where p is some prime number > 3. They are designed in such a way as to maximize the probability that all sequences will converge, while excluding all factors of a given prime number.

The program that excludes factors of 5 appears to have only one loop that has 1 iteration and bottoms out at 3, not 1. It is the exact same as P2 above with a multiply by 5 at an n divisible by 3 and not a multiply by 7, i.e.:

If even: /2; ElseIf divisible by 3: (5n/3 + 1)/2; Else: (5n + 1)/2

A program that looks like the original Collatz but with a 5/2 step instead of a 3/2 step diverges for most n0, as one can imagine since you get an average increase of 5/4 every 2 steps, instead of a decrease of ¾, assuming 5/2 steps are just as likely as 1/2. You can see that since all odd integers are multiplied by 5 before adding the 1, factors of 5 are excluded from the sequences for all n0. In the phrasing of Shaw (2006), all n0 that are factors of 5 will be “pure” in this program. Based on the stability of P2, which converges slower, this program is very stable.

Likewise, we can build a similar program that excludes factors of 7 by adding an extra ElseIf-statement to the program like so:

If even: /2; ElseIf divisible by 3: (7n/3 +1)/2; ElseIf divisible by 5: (7n/5 + 1)/2; Else: (7n+1)/2

This program (that excludes factors of 7) appears to have two loops: the usual 2 iteration loop at n0 = 1 and another 15-iteration loop with a bottom at n0 = 23. Because this program has two loops, it weighs against the idea that excluding a prime factor from all sequences somehow contributes to limiting the program to only one loop. The presence of yet another loop at n0 = 23 (see Table 3) is another mystery.

If we keep expanding the ElseIf statements to create programs that exclude prime factors ≥ 11 we see that they do not converge for most n0. The most simple null model formulation we have been using shows that P6 excluding factors of 7 decreases by an average of only 10% (a factor of 0.9) every 10 steps, and P6 excluding factors ≥ 11 increases with more iterations on average, so they diverge.

**P9**: A Collatz-like closely related set of programs that deviate from Collatz only when n is divisible by selected prime numbers, and only deviate by changing the +/- 1, not the factor of 3. In other words:

If even: /2; ElseIf divisible by p: (3n -/+ 1)/2; Else: (3n +/- 1)/2

The sign of the first -/+ 1 is kept opposite of the sign of the second +/- 1, and they are never set to the same sign (so that no versions of P9 are equivalent to the original Collatz or the 3n – 1 version previously discussed). Thus P9 explores the effect of the +/- 1, not the multiplying factors. Notice in the summary of results in Table 6 that the loops for the 3n – 1 version of Collatz at 1, 5, and 17 are prominent, but new loops at 11 and 125 also appear. The longest loop in Table 6 is L125 with 18 iterations and it rises to a maximum height of 946.

| Prime $p_i$ | L1 | L5 | L11 | L17 | L125 |
|---|---|---|---|---|---|
| **5** | -/+ and +/- | -/+ | | | |
| **7** | -/+ and +/- | | | | |
| **11** | -/+ and +/- | +/- | | +/- | +/- |
| **13** | -/+ and +/- | +/- | -/+ | | |
| **17** | -/+ and +/- | +/- | -/+ | | |
| **19** | -/+ and +/- | +/- | | +/- | |

**Table 6** The table shows where loops exist in different versions of P9, where pi is the prime number used in the first ElseIf of the program, and -/+ refers to the first and second decrement/increment by 1 after the multiply by 3 (-/+ indicates subtract 1 then add 1, +/- is the reverse). There is a loop at L1 in both versions for all $p_i$.

**Q5**: Are there ways of determining, from the definition of a simple program that generates a set of sequences, how they will behave in the 4 qualities that make sequences “interesting” (see **Q4** above), particularly number quality 4? Specifically, can we predict how many “loops” of integers the sequences generated by a given program will contain, how long the loops will be and where they will bottom out?

**Q5.1**: It appears that most of these simple programs, if they have more loops than the simple one that bottoms at n0 = 1, have an even number of loops. Why is that? The exception appears to be m= 67?…Unless there is a loop that has not been registered?

**Q5.2**: Do all sequences for all P4 programs converge regardless of the value of m? Are there some sequences that converge for any value of m?

It certainly looks that for values of m > 73 some n0 do not converge. It is trivial to show that there are an infinitely many n0 for arbitrary values of m that converge using the picket fence numbers (Table 1). Because the Picket Fence numbers not divisible by 5 will become some power of 2 after a (3n+1)/2 step in any of the P4 programs, they and all nodes attached to them will always converge, no matter how large the value of m.

**Conjectures**

Some of the questions give rise to larger conjectures: speculations about Collatz sequences and Collatz-like programs based on what we have observed above.

**C1**: A sort of meta-Collatz conjecture: there appear to be only a finite set of Collatz-like sequence definitions that are interesting, according to how we define "interesting" in **Q4** above. But even if they are finite in number, are there any particularly interesting ones that we have missed here?

**C2**: (Related to **C1** above) Do these Collatz-like sequences and their loops of integers and tree structures have an analog in physical reality, e.g. in quantum mechanics?

Reasons for thinking these simple programs and their associated sequences might be associated with elementary particle physics:

1) They postulate that time is quantized, which is also a hypothesis in quantum physics.
2) The idea of loops of numbers with discreet values (integers) to which all higher values descend is reminiscent of coherent quantum waves, as in for example the de Broglie atom.
3) There seem to be only a limited number of these basic programs that converge, as there are only a limited number of stable, basic particles.
4) The loops discovered might correspond to fundamental energy levels or masses.
5) The loops could correspond to one dimensional objects, like strings, or even zero dimensional objects that are propelled along a one dimensional "track" by some fundamental field or fields which add energy (the multiply steps) or take it away (divide by 2 steps).
6) The integer loops also fit with the concept of the arrow of time, because they form when the program runs in one direction, but not in the reverse.

These are merely wild speculations, but they should be considered generally: i.e., perhaps not these programs, but others somewhat more complex could reproduce behaviour closely analogous to elementary particles – a thought that others have entertained as well (Wolfram 2002).

**Conclusions**

These explorations into Collatz and Collatz-like sequence generators, although not a proof of the Collatz conjecture, serve to give some intuition on why the conjecture is probably true. Specifically, we suspect the Collatz Conjecture is likely true because:

1) Analysis of the Picket Fence numbers (Figure 2 and Table 1) show that Collatz sequences have infinitely many places to exit the sequence and arrive at 1, and that these exit numbers are accessed by various sequences in a regular and predictable way, given the interconnected, tree-like structure of the sequences.
2) Analysis of so-called islands of persistence in Collatz sequences for sets of n0 spaced apart by powers of 2 and 3 – where the sequence lengths are exactly the same for all n0 that differ by hundreds of powers of 2 and 3 – show how densely connected the branches and sub-branches of Collatz sequence trees are. By the time we reach the known lower bound for any possible violation of the Collatz conjecture, the thicket of branches and limbs coming off these islands of persistence must be so dense that it seems impossible that any violation of the conjecture, whether a branch reaching up to infinity or a loop of integers with many 100,000's of numbers, could avoid contact with this dense thicket of branches reaching up to infinity.
3) The closely related P4 programs continue to converge for a wide variety of tested n0, even for values of m well beyond the m = 15 that corresponds to the original Collatz conjecture. Even for m = 73, almost 5x greater than m = 15, all (but one) known sequences converge, although some are incredibly long. This argues that Collatz sequences very likely never diverges for any n0, in keeping with current expectations given the null model.
4) Likewise in the P6 programs, i.e. the programs that exclude prime integer factors from their sequences the way that Collatz sequences exclude factors of 3, and of which the original Collatz program is the fastest converging, also appear to always converge for prime integers 5 and 7 as well, even though these converge much more slowly.
5) In none of the related programs can we find any loops that have a lowest iteration > 100,000, whereas current proofs have shown if Collatz sequences have any other loops than the loop bottoming at 1, they rest orders of magnitude > 100,000.
6) In none of the related programs can we find any loops that have a total length > 1,000 iterations, whereas current proofs have shown if Collatz sequences have any other loops than the one bottoming at 1, they are orders of magnitude longer than 1,000 iterations.
7) Analysis of Collatz-related programs and sequences reveals the fundamental dynamics of these sequences to be a branching process. When the sequence generating programs are run in reverse, we see a branching process that tends to "reach" or incorporate more and more nodes (starting integers) higher up the positive number line of integers. With a finite number of loops at lower integers, the branches apparently become dense enough to capture all integers all the way up the number line, as we discuss in point 2 above.

Much of the language that I use here and that is used in general to discuss and analyze the Collatz conjecture and Collatz-like sequence generators is fundamentally (at least) two dimensional: loops, trees, branches, and graphs. But the sequences themselves and the positive number line that is mapped to them are fundamentally one dimensional. There is a lot of literature in data science devoted to helping people see multi/infinite dimensional processes in lower dimensional spaces like two and three for ease of visualization. Much less, however, is devoted to expanding a one-dimensional process into higher dimensions in ways that faithfully reproduces the constraints of the single dimension. I believe this may be one of the main reasons for our failure to prove the Collatz conjecture to date.

Many of the Collatz-like sequence generating programs manage to shatter the positive number line into several distinct trees which never touch each other. Why the Collatz and several other Collatz-like sequence generators do not is, perhaps, a more interesting question than the Collatz conjecture itself.

**New Queries and Conjectures that Emerge from this Analysis**

Some of these queries and conjectures are discussed above in the alternate Collatz-like sequence generators, but we will summarize and label them here for clarity.

**Query 1** Why do some programs contain only one loop while others contain more?

**Query 2** What is the last m value in program P4 for which all n0 converge to one of several loops? It appears that m > 73, but a proof of the integer limit is lacking.

**Query 3** Why are most sequence lengths emerging from the programs of P4 inferior to even the simplest null model? What allows most P4 sequences to converge relatively fast? Why does P4 continue to converge even as high as m=73? Do all n0 converge for m=73?

**Query 4** How can we predict the number, location, and size of loops in a given program, just from the structure and parameters of that program?

**Query 4.1** Why do most (all?) of the programs that have > 1 loop have an even number of loops?

**Query 4.2** Do all P4 for arbitrary values of m have at least some sequences that converge? Yes, as addressed in the discussion above. There are an infinitely many n0 that converge for all P4, even for values of m where most n0 (probably) diverge.

**Query 4.3** How can we predict the % of positive integers that connect to the nodes for any given loop?

**Conjecture 1** There are only a finite and small number of these simple integer Collatz-like programs that are interesting according to the definition of “interesting” in **Q4** above.

**Conjecture 2** These Collatz-like sequence generating programs are connected in some way to physical reality, e.g. to some quantum mechanical process.

## Methods

Most of the figures, tables, and numerical results in this manuscript were generated using the Mathematica software, version 11.1, on a Windows 10 PC with 16 Gb of RAM. I can send a notebook pertaining to any figure to anyone who requests it. I am indebted to Mathematica StackExchange contributor "mgamer" for an efficient implementation of a basic Collatz sequence generator for Mathematica (https://mathematica.stackexchange.com/questions/69902/collatz-optimization ).

**Disclaimer:** The views expressed in this manuscript are my own, and not to be associated with my employers or professional affiliates. This work was self-financed. Because I am not a professional mathematician, some of the results which I have derived here may be contained elsewhere in the literature without my knowledge. Please forgive my failure to cite those prior results – the oversight was not intentional.

**Acknowledgement:** I would like to thank God and also my family, Jin, Kayla, and Caleb, for their inspiration and for their patience with me.